\documentclass[12pt,twoside]{amsart}  
\usepackage{amssymb,verbatim}  
  
\numberwithin{equation}{section}  
\newtheorem{theorem}{Theorem}[section]  
\newtheorem{lemma}[theorem]{Lemma}  
  
\newtheorem{corollary}[theorem]{Corollary}

\theoremstyle{definition}  
\newtheorem{definition}[theorem]{Definition}
\newtheorem{remark}[theorem]{Remark}

\DeclareMathOperator{\Ann}{Ann}  
\DeclareMathOperator{\Ass}{Ass}  
  
\DeclareMathOperator{\adeg}{adeg}
\DeclareMathOperator{\mult}{mult}

\newcommand{\N}{\mathbb{N}}

\begin{document}  
    
   
\title{Associated primes of monomial ideals and odd holes in graphs}
  
\author{Christopher A. Francisco}  
\address{Department of Mathematics, Oklahoma State University,   
401 Mathematical Sciences, Stillwater, OK 74078}  
\email{chris@math.okstate.edu}  
\urladdr{http://www.math.okstate.edu/$\sim$chris}

\author{Huy T\`ai H\`a}  
\address{Tulane University \\ Department of Mathematics \\  
6823 St. Charles Ave. \\ New Orleans, LA 70118, USA}  
\email{tai@math.tulane.edu}  
\urladdr{http://www.math.tulane.edu/$\sim$tai/}

\author{Adam Van Tuyl}  
\address{Department of Mathematical Sciences \\  
Lakehead University \\  
Thunder Bay, ON P7B 5E1, Canada}  
\email{avantuyl@lakeheadu.ca}  
\urladdr{http://flash.lakeheadu.ca/$\sim$avantuyl/}  
  

\begin{abstract}
Let $G$ be a finite simple graph with edge ideal $I(G)$. Let $I(G)^\vee$ denote the Alexander 
dual of $I(G)$. We show that a description of all induced cycles of odd length in $G$ is 
encoded in the associated primes of $(I(G)^\vee)^2$. This result forms the basis for a method 
to detect odd induced cycles of a graph via ideal operations, e.g., intersections, products 
and colon operations. Moreover, we get a simple algebraic criterion for determining whether a 
graph is perfect. We also show how to determine the existence of odd holes in a graph from the 
value of the arithmetic degree of $(I(G)^\vee)^2$.
\end{abstract}
   
\maketitle  
  
  
\section{Introduction} \label{s.intro}  

A recent breakthrough in graph theory is the Strong Perfect Graph Theorem, 
proved by Chudnovsky, Robertson, Seymour and Thomas \cite{CRST}. This result shows 
that a graph is perfect if and only if neither it nor its complementary graph has an odd 
induced cycle of length at least five; these cycles are often referred to as {\bf odd holes}. 
Consequently, the ability to detect odd induced cycles in a graph in systematic ways is significant. 
Chudnovsky, Cornu\'ejols, Liu, Seymour, and Vu\v{s}kovi\'c \cite{CCLSV} proved the existence of a 
polynomial time algorithm to determine if a graph is perfect.  
However, if $G$ is not perfect, then this algorithm does not 
tell us whether it is $G$ or $G^c$ that contains an odd hole.  More recently, Conforti, 
Cornu\'ejols, Liu, Vu\v{s}kovi\'c, and Zambelli \cite{CCLVZ} showed that one can determine 
if a graph has an odd hole in polynomial time provided $G$ has a bounded clique number.  
In general, there is no known effective algorithm for detecting the existence of odd holes.

Our goal in this paper is to expand the dictionary between graph theory and commutative 
algebra by providing simple, explicit ways to detect all odd induced cycles in graphs, 
allowing us to determine whether a graph is perfect and, if not, where the offending odd 
hole lies. The novelty of our work is the surprising connection between odd holes on the 
graph-theoretic side and associated primes on the commutative algebra side. While the 
algorithms have exponential running time in the worst case, we hope that our results 
will be useful from a theoretical perspective.

More precisely, suppose that $G = (V_G,E_G)$ is a finite simple graph on the vertex set 
$V_G = \{x_1,\ldots,x_n\}$ with edge set $E_G$. By identifying the vertices with the 
variables in the polynomial ring $R =k[x_1,\ldots,x_n]$ over the field $k$, one 
can associate to $G$ a square-free quadratic monomial ideal 
\[I(G) = \left(\{x_ix_j ~|~ \{x_i,x_j\} \in E_G \}\right). \] 
The ideal $I(G)$ is called the {\bf edge ideal} of $G$.  The edge ideal $I(G)$, which was first 
introduced by Villarreal \cite{V2}, is an algebraic object whose invariants can be related to the 
properties of $G$, and vice-versa. Simple graphs and hypergraphs can also be viewed as 
clutters, and so, edge ideals of clutters can be defined in the same way 
(cf. \cite{DV, HM, HMV}). Many researchers have been interested in using the edge ideal 
construction to build a dictionary between the fields of graph theory and commutative 
algebra. For general references, see \cite{SVV,V1,V2}; for invariants encoded in the 
resolution, see \cite{CN,EGHP,HVT,HVT1,HVT2,JK,RVT}; for classes of (sequentially) 
Cohen-Macaulay graphs, see \cite{FH,FV,HHZ,VTV}.

The first main result of this paper is to show that every odd induced cycle in a graph can be detected from
the associated primes of $R/(I(G)^{\vee})^2$, where $I(G)^{\vee}$ is the Alexander dual of $I(G)$, thus 
giving us a method for determining perfection. In fact, 
not only do the associated primes tell us if an odd hole exists, they also 
indicate which vertices make up the odd hole. In particular, we show:

\begin{theorem}[{Corollary \ref{c.primesofsquare}}] \label{maintheorem1}
Let $J = I(G)^{\vee}$. A prime ideal $P=(x_{i_1},\dots,x_{i_s})$ is in $\Ass(R/J^2)$,
the set of associated prime ideals of $R/J^2$, if and only if:  
\begin{enumerate}  
\item $P=(x_{i_1},x_{i_2})$, and $\{x_{i_1},x_{i_2}\}$ is an edge of G, or  
\item $s$ is odd, and the induced graph on $\{x_{i_1},x_{i_2}, \ldots, x_{i_s}\}$ is an induced cycle of $G$.  
\end{enumerate}  
\end{theorem}  

Theorem \ref{maintheorem1} is not the first time the induced odd cycles of a graph 
have been found using commutative algebra.
Simis and Ulrich \cite{SU} showed that $I(G)^{\{2\}}$, the join of $I(G)$ with itself, is generated by the 
square-free monomials $x_{i_1}x_{i_2} \cdots x_{i_r}$, where $r$ is odd, and the induced graph on 
$\{x_{i_1},x_{i_2}, \ldots, x_{i_r}\}$ is an induced cycle.  We find (Theorem \ref{irreducible}) 
an irreducible decomposition for the ideal $(I(G)^{\vee})^2$, and then pair this decomposition with a 
result of Sturmfels and Sullivant \cite{SS} to recover Simis and Ulrich's result (Corollary \ref{2secant}).

Our proof of Theorem \ref{irreducible}, and subsequently, Theorem \ref{maintheorem1}, is 
based upon the notion of a $2$-cover of a graph. If $G$ is simple graph on $n$ vertices, then 
$a = (a_1,\ldots,a_n) \in \N^n$ is a {\bf $2$-cover} if $a_i + a_j \geq 2$ for all edges 
$\{x_i,x_j\} \in E_G$. A $2$-cover $a$ is {\bf reducible} if there exist $b,c \in \N^n$ such that $a = b+c$, 
where $b$ and $c$ are both 1-covers, or one is a 2-cover, and the other is a 0-cover. 
(A 0-cover is simply any nonzero vector $a \in \N^n$.)  Otherwise, we say $a$ is {\bf irreducible}. If we let 
$J = I(G)^{\vee}$, then each generator of  $J^{(2)}$, the second symbolic power of $J$, 
corresponds to some $2$-cover of $G$, while the generators of $J^2$ correspond only to the 
reducible $2$-covers that are the sum of two 1-covers.  The key ingredient in our proof is 
Dupont and Villarreal's 
classification of irreducible $2$-covers \cite{DV}. This classification allows us to describe 
the irreducible decomposition of $J^2$.

The algebra of vertex covers of a (hyper)graph was first studied by Herzog, 
Hibi and Trung \cite{HHT}. In \cite{HHT}, the authors use the terminology of 
{\it decomposable} and {\it indecomposable} covers instead of reducible and irreducible
covers.  We choose to use reducible and 
irreducible covers to be consistent with the result of \cite{DV} that we require.

Theorem \ref{maintheorem1} also forms the basis for other methods to detect the existence of odd holes
in a graph using only the operations of commutative algebra. In particular, the existence of odd holes
can be characterized algebraically as follows:

\begin{theorem}[{Theorems \ref{saturation} and \ref{adeg}}] \label{thm.intro2}
Let $G$ be a simple graph with edge ideal $I(G)$. Set $J = I(G)^{\vee}$,
and let  
\[L = \prod_{1 \leq i_1 < i_2 < i_3 < i_4 \leq n} (x_{i_1}+x_{i_2}+x_{i_3} + x_{i_4}).\]
Then the following are equivalent:
\begin{enumerate}
\item[(a)] $G$ has no odd hole. 
\item[(b)] $J^2:(L) = J^2$.
\item[(c)] $\adeg(J^2) = 3|E_G| + t(G)$, where $\adeg(J^2)$ denotes
the arithmetic degree of $J^2$, and $t(G)$ is the number of triangles of $G$.
\end{enumerate}
\end{theorem}  

\noindent 
The proof of the equivalence of (a) and (b) follows from a well-known lemma  
(see Lemma \ref{decompositionlemma}) 
that the saturation of an ideal $J$ by an ideal $K$ results in an ideal whose
 associated primes do not contain $K$. For the equivalence of (a) and (c), we 
use Sturmfels, Trung, and Vogel's notion of standard pairs 
to compute the arithmetic degree \cite{STV}. 

We note that the main bottleneck in using Theorems ~\ref{maintheorem1} and \ref{thm.intro2}  to determine 
perfection occurs in computing $J$. The generators of $J$ are the minimal vertex 
covers of $G$, and determining the vertex covers of a graph is an NP-complete 
problem \cite{M2book}. Despite this, an algorithm based on Corollary~\ref{c.primesofsquare} 
is simple to code and generally runs reasonably quickly in Macaulay 2. See 
Section~\ref{s.oddcycles} for some specific examples.

The structure of our paper is as follows.  In Section~\ref{s.prelims} we collect together the needed graph theoretic and algebraic results, and we introduce Dupont and Villarreal's classification of irreducible $2$-covers. Section~\ref{s.oddcycles} is devoted to the proof of Theorem \ref{maintheorem1}. 
Finally, Section~\ref{s.otheralgorithms} contains two algebraic characterizations of graphs with odd holes. 
\vspace{.25cm}

\noindent
{\bf Acknowledgments.} The first and third authors began this project at the 
``Syzygies and Hilbert Functions'' workshop at the Banff International Research Station (BIRS), and we 
completed the project during a Research in Teams week at BIRS.  We thank the organizers for inviting 
us to the workshop and BIRS for its hospitality.  The genesis for our results was a number of 
computer experiments using both CoCoA \cite{C} and Macaulay 2 \cite{M2}. 
We also thank Giulio Caviglia, Jeff Mermin, and Susan Morey for some fruitful discussion and Rafael 
Villarreal for sending us an early version of his preprint. Finally, we thank two anonymous referees for their suggestions for improvement. The first author is partially supported by an NSA Young Investigator's Grant and an Oklahoma State University Dean's Incentive Grant. The 
second author is partially supported by Board of Regents Grant LEQSF(2007-10)-RD-A-30 and Tulane's 
Research Enhancement Fund. The third author acknowledges the support provided by NSERC.

  
\section{Graph Theory and Irreducible Covers} \label{s.prelims}
  
In this section we recall the needed terms and results from graph theory, and furthermore, we 
introduce Dupont and Villarreal's classification of irreducible of 2-covers \cite{DV},
which forms a key ingredient for our proof of 
Theorem \ref{irreducible} and Corollary \ref{c.primesofsquare}. 
We continue to use the definitions 
and terms 
from the introduction.
  
Let $G = (V_G,E_G)$ denote a finite simple graph (no loops or multiple edges) on the vertex set 
$V_G = \{x_1,\ldots,x_n\}$ and edge set $E_G$.  We shall abuse notation and write $x_ix_j$ for an 
edge $\{x_i,x_j\} \in E_G$.  If $S \subseteq V_G$, the {\bf induced subgraph} of $G$ on $S$, 
denoted by $G_S$, is the graph with vertex set
$S$ and edge set $E_{G_S} = \{x_ix_j \in E_G ~|~ \{x_i,x_j\} \subseteq S\}$.

\begin{definition} \label{d.cycle} A \textbf{cycle} in a simple graph $G$ is an alternating sequence of 
distinct vertices and edges $C = x_{i_1}e_1x_{i_2}e_2 \cdots x_{i_{n-1}}e_{n-1}x_{i_n}e_nx_{i_1}$ in 
which the edge $e_j$ connects the vertices $x_{i_j}$ and $x_{i_{j+1}}$ ($x_{i_{n+1}} = x_{i_1}$) 
for all $j$.   In this case, $C$ has \textbf{length} $n$ and we call $C$ an $n$-cycle. 
We shall often write a cycle simply as $x_{i_1}x_{i_2} \cdots x_{i_n}x_{i_1}$ or $x_{i_1} \cdots x_{i_n}$, omitting the edges. A {\bf chord} is an edge that joins two nonadjacent vertices in the cycle.  We shall 
use $C_n$ to denote an $n$-cycle without any chords.  We usually
refer to $C_n$ as an {\bf induced cycle} since the induced graph on $\{x_{i_1},x_{i_2}, \ldots, x_{i_n}\}$
contains only the edges and vertices in the cycle.
If an induced cycle has an odd (resp. even) number of 
vertices, we shall call it an {\bf odd} (resp., {\bf even}) {\bf cycle}.  An odd induced cycle of length at least 
five is called an {\bf odd hole}.
\end{definition}  

A subset $W$ of $V_G$ is a {\bf vertex cover} of $G$ if every edge is incident 
to at least one vertex of $W$.  A vertex cover $W$  
is a {\bf minimal vertex cover} if no proper subset of $W$ is a vertex cover.  
More generally, we can define vertex covers of any order.
  
\begin{definition} Let $\N$ denotes the set of nonnegative integers. If $G$ is a simple graph on $n$ vertices, then a nonzero vector $a = (a_1,\ldots,a_n) \in \N^n$ is a 
{\bf vertex cover of order $k$} (or a {\bf $k$-cover}) if $a_i + a_j \geq k$ for all edges $x_ix_j \in E_G$.  
A $k$-cover $a$ is {\bf reducible} if there exists an $i$-cover $b \in \N^n$ and a $j$-cover $c \in \N^n$  
such that $a = b+c$ and $k = i+j$.  Otherwise, we say $a$ is {\bf irreducible}.  
\end{definition}  
  
\begin{remark} When $(a_1,\ldots,a_n)$ is a $(0,1)$-tuple,  
then a vertex cover of order $1$ corresponds to the standard notion  
of a vertex cover.  
At times, we shall write vertex covers of order $k$ as monomials, using the usual 
correspondence between monomials and vectors of nonnegative integers, i.e.,
$x_1^{a_1} \cdots x_n^{a_n}$ corresponds to the cover $(a_1, \dots, a_n)$.
\end{remark}

We shall use a result of Dupont and Villarreal \cite{DV} that classifies irreducible $2$-covers.  
In fact, we shall state a slightly more general version than their result. 
The proof of the theorem was already embedded in the proof of \cite[Theorem 5.1]{HHT}. 
In the statement below, the set $N(A)$ denotes the {\bf neighbors} of the
set $A \subseteq V_G$, that is, \[N(A) = \{y \in V_G\setminus A ~|~ \mbox{there exists $x \in A$
such that $xy \in E_G$}\}.\]
A set of vertices $A \subseteq V_G$ is {\bf independent} if the induced
graph $G_A$ contains no edges, that is, there are no edges among
the vertices of $A$.
As well, $G$ is called {\bf bipartite} if we can partition $V_G = V_1 \cup V_2$
so that every $xy \in E_G$ has the property that $x \in V_1$ and $y \in V_2$. 

\begin{theorem}[see {\cite[Theorem 2.6]{DV}}]\label{irreducible2}  Let $G$ be a simple graph.  
\begin{enumerate}  
\item[$(i)$]  If $G$ is bipartite, then $G$ has no irreducible $2$-covers.  
\item[$(ii)$]  If $G$ is not bipartite and $a$ is a $2$-cover that cannot be written as the sum of two $1$-covers, then (up to some permutation of the vertices)  
\[a = (\underbrace{0,\ldots,0}_{|A|},\underbrace{b_1,\ldots,b_{|B|}}_{|B|},1,\ldots,1)\]  
for some (possibly empty) independent set $A$ and a set $B \supseteq N(A)$ such that   
\begin{enumerate}  
\item[(1)] $b_j \ge 2$ for all $j = 1, \dots, |B|$,
\item[(2)] $B$ is not a vertex cover of $G$ and $V \neq A \cup B$, and  
\item[(3)] the induced subgraph on $C = V \setminus (A \cup B)$ is not bipartite. 
\end{enumerate}  
Moreover, if $a$ is irreducible, then 
\[a = (\underbrace{0,\ldots,0}_{|A|},\underbrace{2,\ldots,2}_{|B|},1,\ldots,1),\]
$B = N(A)$, and the induced subgraph on $C$ has no isolated vertices.
\end{enumerate}   
\end{theorem}  

\begin{proof} Part $(i)$ follows from part (b) of \cite[Theorem 5.1]{HHT}. To prove part $(ii)$, we let 
$A$ be the set of vertices $x_i$ such that $a_i = 0$. Since $a$ is a $2$-cover, for any $x_i \in N(A)$,
 we must have $a_i \ge 2$. We may also include in $B$ all other vertices $x_j$ not in $N(A)$ such that 
$a_j \ge 2$. Clearly, $B \supseteq N(A)$, and (1) is satisfied.

If $B$ is a vertex cover, then $(0,\ldots,0,c_1,\ldots,c_{|B|},d_1,\ldots,d_{|C|})$, where $c_i \ge 1$ 
and $d_j \ge 0$ for all $i$ and $j$, is a $1$-cover of $G$. Thus, $a$ can be written as the sum of two 
$1$-covers, a contradiction. Therefore, $B$ is not a vertex cover of $G$. This also implies that $C$ is
 not empty, and (2) is satisfied.

It follows from part (b) of \cite[Theorem 5.1]{HHT} that if the induced subgraph on $C$ is bipartite, 
then it admits $(1,\dots,1)$ as a $2$-cover that can be written as the sum of two $1$-covers. Moreover, 
$(\underbrace{0,\ldots,0}_{|A|},\underbrace{c_1,\ldots,c_{|B|}}_{|B|})$, where $c_i \ge 1$ for all $i$,
 is a $1$-cover of the induced subgraph on $A \cup B$. Thus, $a$ can be written as the sum of two 
$1$-covers, a contradiction. Thus (3) is satisfied.

To prove the last statement we observe that if $b_j > 2$ for some $j$, then $a$ can be written as the 
sum of a $0$-cover 
$(\underbrace{0,\ldots,0}_{|A|},\underbrace{0, \dots, 0, 1, 0 \dots, 0}_{1 \text{ at the $j$-th place }},\underbrace{0,\ldots,0}_{|C|})$ 
and another $2$-cover. This contradicts the irreducibility of $a$. 
Also, if there exists some $x_j \in B \setminus N(A)$, then $a$  can be written
in a similar fashion as the sum of a $0$-cover and a $2$-cover.  This contradiction
thus implies that $B = N(A)$.
Similarly, if the induced subgraph
 on $C$ has an isolated vertex, say the last one, then $a$ can be written as the sum of the $0$-cover 
$(\underbrace{0,\ldots,0}_{|A|},\underbrace{0,\ldots,0}_{|B|},\underbrace{0,\ldots,0, 1}_{|C|})$ and 
another $2$-cover, again a contradiction.
\end{proof}

We round out this section by explaining how the $2$-covers of a graph $G$ are
related to the Alexander dual of the edge ideal $I(G)$.  First, we define the Alexander dual:

\begin{definition} \label{d.alexdual} Suppose $I$ is a square-free monomial ideal. 
The {\bf Alexander dual} of $I$, denoted by $I^{\vee}$, is the ideal whose primary components are given by the minimal generators of $I$. That is, if $I=(x_{1,1} \cdots x_{1,t_1}, \dots, x_{r,1} \cdots x_{r,t_r})$,
then \[ I^{\vee} = (x_{1,1}, \dots, x_{1,t_1}) \cap \cdots \cap (x_{r,1}, \dots, x_{r,t_r}).\]
\end{definition}  

The ideal $I(G)^{\vee}$ is sometimes referred to as the {\bf cover ideal} because 
of the well-known fact that the generators of $I(G)^{\vee}$ correspond to vertex covers (see, e.g., \cite{HHT}).

Next, we recall the notion of the symbolic power of an ideal, restricting to the case in which 
$I \subset R$ is a square-free monomial ideal. Suppose $I$ has the primary decomposition 
\[ I = P_1 \cap \cdots \cap P_r,\] where each $P_i$ is an ideal generated by a subset of the variables of $R$. 
The {\boldmath $j$}{\bf -th symbolic power} of $I$ is the ideal \[ I^{(j)} = P_1^j \cap \cdots \cap P_r^j.\]

Set $J = I(G)^{\vee}$.  Graph-theoretically, we can interpret the minimal generators of  $J^{(2)}$ and $J^2$ in 
terms of $2$-covers. For convenience, we denote covers by their corresponding monomials instead of the 
vectors themselves. Note that \[ J^{(2)} = \bigcap_{x_ix_j \in E_G} (x_i,x_j)^2,\] so $J^{(2)}$ is the 
ideal whose minimal 
generators yield a 2-cover of $G$. On the other hand, $J^2$ is more restrictive. Its minimal generators
 are still 2-covers, but they must be able to be partitioned into two ordinary vertex covers. That is, 
if $m \in J^2$, 
then $m=m'm''$, where $m'$ and $m''$ are 1-covers of $G$. A main part of the proof of 
Theorem \ref{irreducible} is to understand, via Theorem \ref{irreducible2}, the difference between 
the monomials in $J^{(2)}$ and those in $J^2$.

  
\section{Odd cycles and associated primes} \label{s.oddcycles}

In this section we prove the main result of our paper, that is, the odd cycle structure of a graph $G$ appears 
in the associated primes of $R/J^2$, where $J=I(G)^{\vee}$.

The following definition is classical:

\begin{definition} \label{d.assocprime} Let $M$ be an $R$-module. A prime ideal $P$ is called an 
{\bf associated prime} of $M$ if $P=\Ann(m)$, the annihilator of $m$, for some $m \in M$. The set of 
all associated primes of $M$ is denoted by $\Ass(M)$.
\end{definition}  

We begin with some observations.  Because we will only be dealing with the case in which $I$ is a
monomial ideal, all $P \in \Ass(R/I)$ will have the form $P = (x_{i_1},\ldots,x_{i_t})$ for some 
subset $\{x_{i_1},\ldots,x_{i_t}\} \subset \{x_1,\ldots,x_n\}$. Since 
$J = I(G)^{\vee} = \bigcap_{x_ix_j \in E_G} (x_i,x_j),$  
the associated primes of $R/J$ are exactly the primes corresponding to the edges of $G$, that is, 
the prime ideals $(x_i,x_j)$ where $x_ix_j$ is an edge of $G$.  Moreover, $\Ass(R/J^{(2)})=\Ass(R/J)$. 
However, $R/J^2$ can have additional associated primes, and it is these primes we seek to identify.  
We proceed by computing something stronger, namely, an irreducible decomposition for $J^2$. An 
{\bf irreducible} monomial ideal in $n$ variables is an ideal of the form $(x_1^{a_1},\dots,x_n^{a_n})$ 
with ${\bf a} = (a_1,\dots,a_n) \in \N^n$. This ideal is usually denoted as ${\bf m}^{\bf a}$, so, 
for example, the maximal homogeneous ideal would be ${\bf m}^{(1,\dots,1)}$. If $a_i=0$, then we adopt the convention that ${\bf m}^{\bf a}=(x_1^{a_1}, \dots, \widehat{x_i^{a_i}}, \dots, x_n^{a_n})$; that is, no power of $x_i$ is in the ideal. Every monomial ideal $I$ 
can be decomposed into the intersection of finitely many irreducible ideals, 
i.e., $I = {\bf m}^{{\bf a}_1} \cap \cdots \cap {\bf m}^{{\bf a}_s}$ (see, for example, \cite[Lemma 5.18]{MS}).

\begin{theorem}\label{irreducible}
Let $G$ be a finite simple graph.  If $J = I(G)^{\vee}$, then the irredundant irreducible decomposition of $J^2$ is:
\[J^2 = \bigcap_{x_ix_j \in E_G} [(x_i^2, x_j) \cap (x_i,x_j^2)]  \cap \bigcap_{\begin{tabular}{c}
\mbox{induced graph on $\{x_{i_1},\ldots,x_{i_s}\}$} \\
\mbox{ is an odd cycle}\end{tabular}} (x_{i_1}^2,\ldots,x_{i_s}^2).\]
\end{theorem}

\begin{proof}
Let $L$ denote the ideal on the right-hand side in the statement of the theorem.
  
Consider a minimal generator  $M \in J^2$, and thus $M$ is the product of two 
1-covers of $G$.  
Since $I(G)^{\vee} \subseteq (x_i,x_j)$ for every $x_ix_j \in E_G$, it follows that 
$J^2 \subseteq (x_i,x_j)^2 \subseteq (x_i^2,x_j)$,
and similarly, $J^2 \subseteq (x_i,x_j^2)$. Hence 
\[M \in \bigcap_{x_ix_j \in E_G} [(x_i^2,x_j) \cap (x_i,x_j^2)].\]

Suppose that the graph $G$ has an odd induced cycle on the vertices 
$\{x_{i_1},\ldots,x_{i_s}\}$.  We claim
that there exists some $x_{i_j} \in \{x_{i_1},\ldots,x_{i_s}\}$ such that 
$x_{i_j}^2\mid M$.  Suppose not.  Then  $M = x_{i_1}^{a_{i_1}}\cdots x_{i_s}^{a_{i_s}}M'$, where 
$0 \leq a_{i_j} \leq 1$ for 
all $j = 1,\ldots,s$,  and
no $x_{i_j}$ divides $M'$.  Since $M = M_1M_2$, where $M_1,M_2 \in J$, both $M_1$ and $M_2$ must 
contain at least $(s+1)/2$ vertices of $\{x_{i_1},\ldots,x_{i_s}\}$ in order to 
cover the odd induced 
cycle on these vertices.  So, in the variables $\{x_{i_1},\ldots,x_{i_s}\}$, $M$
 must have degree at 
least $s+1$.  But we have assumed that
$M$ has degree at most $s$ in these variables, a contradiction.  So, there exists some $x_{i_j}$ 
such that $x_{i_j}^2\mid M$.  Hence $M \in (x_{i_1}^2,\ldots,x_{i_s}^2)$.  Because this is true for 
each odd induced cycle, $M$ is also in the second set of intersections. Thus, 
$J^2 \subseteq L$.

We now prove the converse. Consider a minimal generator $N$ of $L$. Since 
$N \in \bigcap_{x_ix_j \in E_G} [(x_i^2,x_j) \cap (x_i,x_j^2)]$, it is clear that 
$N$ is a $2$-cover of $G$. 
It suffices to show that $N$ can be written as the sum of two $1$-covers. Suppose that this is not the case. 
Then by Theorem \ref{irreducible2}, for some independent set 
$A$, $B \supseteq N(A)$ and $C = V \setminus (A \cup B) \not= \emptyset$, we have 
\[N = \prod_{x_j \in B} x_j^{b_j} \prod_{x_j \in C} x_j\] where $b_j \ge 2$ for 
all $j$, and the 
induced subgraph on $C$ is not a bipartite graph. This 
implies that the induced subgraph on $C$ contains an odd cycle, say on the vertices $\{x_{i_1}, \dots, x_{i_s}\}$. 
From the expression of $L$, we have $N \in (x_{i_1}^2, \dots, x_{i_s}^2)$. However, as we have seen, the 
power of any vertices of $C$ in $N$ is exactly one. This is a contradiction. Hence, any
minimal generator of $L$ can be written as 
the sum of two $1$-covers, that is, $L \subseteq J^2$. 
\end{proof}

\begin{remark} \label{r.higherpowers}
The irreducible decomposition of $J^s$, which is studied in \cite{FHVT}, 
is substantially more complicated when $s > 2$ because of its 
relation to the chromatic number of 
induced subgraphs. When $G$ is an odd cycle, the situation is much easier; the associated 
primes of $J^s$ correspond to the edges and the odd cycle itself, and the exponents that appear 
in the irreducible components come from a formula depending on $s$
and $|V_G|$.
\end{remark} 
 
Our main result is now an immediate corollary of Theorem~\ref{irreducible}.

\begin{corollary} \label{c.primesofsquare}  
Let $G$ be a finite simple graph, and set $J = I(G)^{\vee}$.
A prime $P=(x_{i_1},\dots,x_{i_s})$ is in $\Ass(R/J^2)$ if and only if:  
\begin{enumerate}  
\item $P=(x_{i_1},x_{i_2})$, and $x_{i_1}x_{i_2}$ is an edge of G, or  
\item $s$ is odd, and after re-indexing,   
$x_{i_1}x_{i_2} \cdots x_{i_s}x_{i_1}$ is an induced cycle of $G$.  
\end{enumerate}  
\end{corollary}    

Moreover, we get a method for detecting perfect graphs from the following corollary:

\begin{corollary}\label{c.perfect}
Let $G$ be a finite simple graph with $J=I(G)^{\vee}$ and $J_c=I(G^c)^{\vee}$, where $G^c$ is the 
complementary graph of $G$. Then $G$ is perfect if and only if neither $\Ass(R/J^2)$ nor $\Ass(R/J_c^2)$ 
contains a prime of height greater than three.
\end{corollary}

Our intent in this paper is to focus primarily on the connections between commutative algebra and 
graph theory and not on the speed of algorithms. However, we make a few comments here about
using Corollary~\ref{c.perfect} to detect perfect graphs.   While any algorithm based on
Corollary \ref{c.perfect} does not run in polynomial time, it 
has the advantage that it tells us exactly where any odd holes occur, and whether they are in 
$G$ or $G^c$, which the polynomial time algorithm from \cite{CCLSV} does not.
An algorithm based upon Corollary ~\ref{c.perfect} has been included
in the Macaulay 2 package EdgeIdeals \cite{EdgeIdeals} under the command {\tt allOddHoles}.
 Moreover, despite 
bad worst-case running time, {\tt allOddHoles} has been successful in 
computing relatively large examples. We ran three examples on a standardly-equipped laptop
 and computed the time this algorithm took to detect all the odd 
holes in randomly-chosen graphs. For a randomly-chosen graph on 14 vertices with 40 edges, 
{\tt allOddHoles} took 0.047 seconds to identify all odd holes in the graph. The command took 
0.468 seconds for a randomly-chosen graph on 20 vertices with 60 edges, and it took 13.665 
seconds for a randomly-chosen graph on 30 vertices with 200 edges. These times strike us as 
reasonable given the size of the graphs and the usual difficulties of working with large polynomial 
rings in computer algebra systems. 

Corollary \ref{c.primesofsquare} also provides some crude bounds on the depth and projective dimension
of $R/J^2$ in terms of the size of the largest induced odd cycle.

\begin{corollary}  Let $G$ be a finite simple graph on $n$ vertices, and let $t$ denote
the size of the largest induced odd cycle of $G$.  If $J = I(G)^{\vee}$, then
\begin{enumerate}
\item[(a)] $\operatorname{depth}(R/J^2) \leq n -t$,
\item[(b)] $\operatorname{projdim}(R/J^2) \geq t$.
\end{enumerate}
\end{corollary}

\begin{proof} By the Auslander-Buchsbaum Formula, it suffices to prove (a).  
For any ideal $I$ of $R$, $\operatorname{depth}(R/I) \leq \dim{R/P}$ for any $P \in \Ass(R/I)$.
Now apply Corollary \ref{c.primesofsquare}.
\end{proof}

We round out this section by using our methods to give an alternate proof of a 
result of Simis and Ulrich \cite{SU} 
and Sturmfels and Sullivant \cite{SS} about the second-secant of $I(G)$.

The join of an ideal was studied in \cite{SU} and \cite{SS}.  We recall a special case of this definition.
If $I$ and $J$ are ideals of $k[x_1,\ldots,x_n]$, then their {\bf join}, denoted
$I * J$, is a new ideal of $k[x_1,\ldots,x_n]$ which is computed as follows:
Introduce new variables $y = \{y_1,\ldots,y_n\}$ and $z = \{z_1, \dots, z_n\}$, and let $I(y)$ (resp. $J(z)$)
denote the image of the ideal $I$ (resp. $J$) under the map $x_i \mapsto y_i$ (resp. $x_i \mapsto z_i$) in
the ring $k[x_1, \dots, x_n, y_1, \dots, y_n, z_1, \dots, z_n]$.  Then 
\[I * J = (I(y)+J(z) + (y_1+z_1-x_1, \dots,y_n+z_n -x_n)) \cap k[x_1,\ldots,x_n].\]
When $I = J$, we call $I * I$ the {\bf second-secant ideal} of $I$ and denote
it by $I^{\{2\}}$.  In the proof of the following theorem,
we use the notation of $I^{[{\bf a}]}$ found in Section 5.2 of \cite{MS}
for the generalized Alexander dual of $I$. 

\begin{corollary}[{\cite[Proposition 5.1]{SU},\cite[Corollary 3.3]{SS}}]\label{2secant}
Let $G$ be a finite simple graph. Then 
\[I(G)^{\{2\}} = (\{ x_{i_1}\cdots x_{i_s} ~|~ \mbox{$G_{\{x_{i_1},\ldots,x_{i_s}\}}$ is an odd induced cycle}\});\]
that is, the generators correspond to the vertices of the induced odd cycles.
\end{corollary}

\begin{proof}
Since $J = I(G)^{\vee}$ is a square-free monomial ideal, every monomial generator of 
$J$ divides $x_1\cdots x_n$, and so every generator of $J^2$ divides $x_1^2\cdots x_n^2$.
  Set ${\bf 1} = (1,\ldots,1) \in \N^n$.  By applying
\cite[Corollary 2.7]{SS} with ${\bf a} ={\bf 1}$, which is sufficiently large since $J$ is square-free, we have
\[I(G)^{\{2\}} = \left((I(G)^{[\bf 1]})^2\right)^{[2\cdot{\bf 1}]} ~~\mbox{modulo}~~ {\bf m}^{{\bf 1}+{\bf 1}},\]
where modulo ${\bf m}^{{\bf 1} + {\bf 1}} = {\bf m}^{{\bf 2}} = (x_1^2,\ldots,x_n^2)$
refers to removing all the monomial generators divisible by $x_i^2$ for some $i$.
Now $I(G)^{[\bf 1]} = I(G)^{\vee}$, so $I(G)^{\{2\}} = (J^2)^{[{\bf 2}]} ~~\mbox{modulo}~~ {\bf m}^{{\bf 2}}$,
where ${\bf 2}: = 2\cdot {\bf 1} = (2,\ldots,2)$.  By \cite[Theorem 5.27]{MS},
the generators of $(J^2)^{[{\bf 2}]}$ are in one-to-one correspondence
with the irreducible components of $J^2$;  in particular, by Theorem \ref{irreducible},
combined with \cite[Theorem 5.27]{MS}, we have
\[(J^2)^{\bf [2]} = (\{x_ix_j^2,x_i^2x_j ~|~ x_ix_j \in E_G\}) + 
(\{x_{i_1}\cdots x_{i_s} ~|~\mbox{$G_{\{x_{i_1},\ldots,x_{i_s}\}}$ is an odd induced cycle}\}).\]
When we remove the monomial generators of $(J^2)^{\bf [2]}$ divisible by $x_i^2$ for some $i$,
we are removing the first ideal, while the second remains, and hence the conclusion follows.
\end{proof}

  
\section{Algebraic classification of odd cycles} \label{s.otheralgorithms}
 
In this section we describe two algebraic approaches to detecting the existence of odd induced cycles 
(and, in particular, odd holes) in a graph. The first method is based upon taking quotients of ideals and is 
well-suited for constructing an algorithm to detect odd cycles
using the ideal operations of commutative algebra.  The second method
is based upon the arithmetic degree of an ideal, which, although
hard to compute, is interesting from a theoretical point of view.
Of course, one could use Corollary \ref{c.primesofsquare} to determine
if a graph has an odd cycle;  however, Corollary \ref{c.primesofsquare}
not only tells us if an odd cycle exists, it tells us
which vertices make up the cycle.  If one is simply
interested in the question of existence, the results of 
this section may be more relevant.

\subsection{Method 1: Colon ideals} Using the technique of ideal saturation, we can describe an algebraic 
approach to detecting odd cycles. Recall that if $I$ and $K$ are ideals of $R$, then the {\bf saturation} 
of $I$ with respect to $K$, denoted $(I:K^{\infty})$, is defined by: 
\[(I:K^{\infty}) =   \bigcup_{N \ge 1} (I:K^N).\]
The ideal $I:K^{\infty}$ is then related to the primary decomposition of $I$ as in 
Lemma \ref{decompositionlemma}. We omit the proof; see, for example \cite[Lemma 2.4]{EHV}.

\begin{lemma}  \label{decompositionlemma}  
Let $I$ be an ideal of $R = k[x_1,\ldots,x_n]$  
with primary decomposition  
\[I = Q_1 \cap Q_2 \cap \cdots \cap Q_r.\]  
If $K$ is an ideal of $R$, then   
\[(I:K^{\infty}) = \bigcap_{K \not\subset \sqrt{Q_i}} Q_i.\]  
\end{lemma}  

We use this result to give a method for detecting odd induced cycles.
  
\begin{theorem}  \label{saturation}  
Let $G$ be a simple graph, and set $J = I(G)^{\vee}$.  Fix an integer $t > 1$, and set 
\[L_t = \prod_{1 \leq i_1 < i_2 < \cdots < i_t \leq n} (x_{i_1}+x_{i_2}+\cdots+x_{i_t}).\]
Then $G$ has no odd induced cycle of length $\geq t$ if and only if $J^2:(L_t) = J^2$.
\end{theorem}  
    
\begin{proof}  Let $J^2 = Q_1 \cap \cdots \cap Q_r$ be the primary decomposition of $J^2$.  By 
Corollary \ref{c.primesofsquare}, we know that $\sqrt{Q_i} = (x_{i_1},x_{i_2})$ where 
$\{x_{i_1},x_{i_2}\}$ is an edge of our graph, or $\sqrt{Q_i} = (x_{i_1},\ldots,x_{i_s})$ with $s$ odd, 
and the induced graph on the vertices in $\sqrt{Q_i}$ is a cycle of odd length. Note that if $K$ and $L$ are any ideals, $K:L^{\infty}=K$ if and only if $K:L=K$, so it suffices to show that $G$ has no odd induced cycle of length $\geq t$ if and only if  $J^2:(L_t)^{\infty} = J^2$.
  
Suppose that $G$ has no odd induced cycle of length $\geq t$, i.e., if $P_i = (x_{i_1},\ldots,x_{i_s})$ 
is an associated prime, then $s=2$ or $s$ is odd and $s < t$.   In both cases 
$(L_t) \not\subset \sqrt{Q_i}= P_i$ for all $i$. Hence, by  Lemma 
\ref{decompositionlemma} 
\[J^2:(L_t)^{\infty} = \bigcap_{(L_t) \not\subset \sqrt{Q_i}} Q_i = \bigcap^r_{i=1} Q_i = J^2.\]  
  
On the other hand, suppose that $G$ has an odd induced cycle of length $\geq t$, i.e., there exists some 
$Q_i$ such that $\sqrt{Q_i} = (x_{i_1},\ldots,x_{i_s})$ with $s \geq t$ odd. Now, 
$x_{i_1}+x_{i_2}+\cdots+x_{i_t} \in \sqrt{Q_i}$, so $L_t \in \sqrt{Q_i}$, and hence 
\[ J^2:(L_t)^{\infty} = \bigcap_{(L_t) \not\subset \sqrt{Q_i}} Q_i \supsetneq  
\bigcap_{i=1}^r Q_i = J^2.\] The result now follows.
\end{proof}  
  
By specializing Theorem \ref{saturation} to the case that $t=4$, we can detect graphs with odd holes:

\begin{corollary}\label{detectoddholes}
Let $G$ be a simple graph, and set $J = I(G)^{\vee}$.  Set 
\[L = \prod_{1 \leq i_1 < i_2 < i_3 < i_4 \leq n} (x_{i_1}+x_{i_2}+x_{i_3}+x_{i_4}).\]
Then $G$ has an odd hole if and only if $J^2:(L) \supsetneq J^2$.
\end{corollary}
  
\subsection{Method 2:  Arithmetic degree}  The second main result of this section is to 
show that one can identify graphs with odd holes via the arithmetic degree.
  
\begin{definition}  
Let $I$ be a homogeneous ideal of $R = k[x_1,\ldots,x_n]$.  The  
{\it arithmetic degree} of $I$ is   
\[\operatorname{adeg}(I) = \sum_{\mbox{homogeneous prime ideals $P \subseteq R$}} 
\operatorname{mult}_I(P)\deg(P).  
\]  
\end{definition}  
  
In the above definition, $\mult_I(P)$ is the length of the largest ideal of finite length  
in the ring $R_P/IR_P$.  It can be shown that $\mult_I(P) > 0$ if and only if $P$ is  
an associated prime of $I$.  So, the above formula gives us information  
about the existence of certain associated primes.  Note that when $I$ is a monomial   
ideal, all the associated primes have the form $P = (x_{i_1},\ldots,x_{i_s})$,  
and $\deg(P) =1$ for all of these ideals.  So, when $I$ is a monomial ideal,  
the above formula reduces to  
\[\adeg(I) = \sum_{\mbox{$P \in \Ass(R/I)$}} \mult_I(P).\]  
  
In the paper of Sturmfels, Trung, and Vogel \cite{STV}, a combinatorial formula for $\mult_I(P)$ is 
given when $I$ is a monomial ideal. Let $X = \{x_1,\ldots,x_n\}$.  Any prime monomial ideal of $R$ 
is generated by some subset of the variables.  In particular, any monomial prime is determined by 
the variables not in the ideal; that is, for each monomial prime ideal $P$, there is a subset $Z \subseteq X$
 such that $P=P_Z := (\{x_i ~|~ x_i \in X \setminus Z\})$.  For a monomial $M \in R$, we let 
$\operatorname{supp}(M)$ denote the {\bf support} of $M$, i.e., the set of variables appearing in
$M$. By \cite[Lemma 3.3]{STV}, $\mult_I(P_Z)$ equals the number of standard pairs  
of the form $(\cdot,Z)$.  If $M$ is a monomial, and $Z \subseteq X$, a pair $(M,Z)$ is {\bf standard} if   
\begin{enumerate}  
\item[(a)]  $(M,Z)$ is {\bf admissible}, i.e., $\operatorname{supp}(M) \cap Z = \emptyset$,  
\item[(b)]  $\left(M\cdot k[Z]\right) \cap I = \emptyset$, and  
\item[(c)] $(M,Z)$ is minimal with respect to the partial order  
\[(M,Z) \leq (M',Z') \Leftrightarrow \mbox{ $M$ divides $M'$ and  
$\operatorname{supp}(M'/M) \cup Z' \subseteq Z$}\]  
for all pairs $(M,Z)$ that satisfy (b).  
\end{enumerate}  
  
We now specialize to the case of the monomial ideal $J^2 = (I(G)^{\vee})^2$.  
By Corollary \ref{c.primesofsquare}, we know that $P_Z$  
is an associated prime of $J^2$ if the induced graph on $X\setminus Z$ is either an edge of  
$G$ or an odd cycle of $G$. Note that if $x_ix_j \in E_G$, and $Z = X\setminus \{x_i,x_j\}$, then 
\[ \mult_{J^2}(P_Z) = \mult_{J^2}((x_i,x_j)) = \deg (x_i,x_j)^2 =3\] 
because $P_Z$ is a minimal prime. (See, e.g., \cite[Section 1]{STV}. We thank an anonymous referee for pointing this out and greatly simplifying our original argument.)

We need one other multiplicity calculation.  

\begin{lemma}  \label{lemma2}
Let the induced graph on $\{x_i,x_j,x_k\}$ be a three-cycle, and set $Z = X \setminus \{x_i,x_j,x_k\}$.  Then 
$\mult_{J^2}(P_Z) = 1.$  
\end{lemma}  
  
\begin{proof}  We begin with some observations about the generators of $J = I(G)^{\vee}$  
and $J^2$.  If $M$ is a minimal generator of $J$, then $M$ must be divisible  
by one of $x_ix_j$, $x_ix_k$ or $x_jx_k$ because $M$ corresponds to a minimal vertex cover,  
and we need at least two of the three vertices to cover the edges of the  
triangle formed by $\{x_i,x_j,x_k\}$.  Hence, any monomial of $J^2$ must be divisible  
by one of
$x_i^2x_j^2,x_i^2x_jx_k,x_ix_j^2x_k,x_i^2x_k^2,x_ix_jx_k^2,x_j^2x_k^2$.  
In particular, every monomial of $J^2$ must be divisible by at least one  
of $x_i^2,x_j^2,x_k^2$.  
  
All the admissible pairs of the form $(\cdot,Z)$ are:   
\[(1,Z),(x_i^a,Z),(x_j^b,Z),(x_k^c,Z),(x_i^ax_j^b,Z),(x_i^ax_k^c,Z),(x_j^bx_k^c,Z),(x_i^ax_j^bx_k^c,Z).\]  
We claim that all but the last pair fail to be a standard pair, and the last is standard only when $a=b=c=1$. 
The conclusion of the lemma then follows.

Note that $(1,Z \cup \{x_i\}) < (1,Z)$, and $(1,Z\cup\{x_i\})$ has the property that  
$k[Z\cup\{x_i\}] \cap J^2 = \emptyset$ since every monomial of $J^2$ is divisible  
by either $x_j$ or $x_k$ (because $x_jx_k$ is an edge of $G$), but no  
such monomial belongs to $k[Z \cup \{x_i\}]$.  So, $(1,Z)$ is not a standard pair since
it is not minimal with respect to the partial order.
  
For $(x_i^a,Z)$, we have $(1,Z\cup\{x_i\}) < (x_i^a,Z)$ since $1|x_i^a$ and 
$\operatorname{supp}(x_i^a/1)\cup Z \subseteq Z \cup \{x_i\}$. But as noted above, 
$k[Z \cup \{x_i\}] \cap J^2 = \emptyset$. Thus $(x_i^a,Z)$ is not minimal, so it cannot be a 
standard pair. A similar argument eliminates $(x_j^b,Z)$ and $(x_k^c,Z)$.  
 
To rule out $(x_i^ax_j^b,Z)$, we first note that $(x_j,Z \cup \{x_i\}) < (x_i^ax_j,Z)$. But we 
also have $x_jk[Z\cup\{x_i\}] \cap J^2 = \emptyset$ since for every monomial $M$ in $J^2$ 
such that $x_j|M$ but $x_j^2 \nmid M$, we must have $x_k|M$. But $x_k \not\in k[Z\cup\{x_i\}]$, 
and hence the intersection is empty. Therefore $(x_i^ax_j,Z)$ is not standard, and by a symmetric 
argument, neither is $(x_ix_j^b,Z)$. Suppose now that $a,b>1$, and consider $(x_i^ax_j^b,Z)$. 
Pick a minimal vertex cover $M_1$ of $G$ containing $x_i$ and $x_j$ but not $x_k$. Then 
$M_1^2=x_i^2x_j^2N_1^2$, where $N_1 \in k[Z]$. Thus 
$x_i^{a-2}x_j^{b-2}M_1^2 = x_i^ax_j^bN_1^2 \in x_i^ax_j^bk[Z] \cap J^2$, and $(x_i^ax_j^b,Z)$ 
is not standard. The same argument with the variables permuted eliminates $(x_i^ax_k^c,Z)$ and $(x_j^bx_k^c,Z)$.

Suppose now that $a>1$, and consider the pair $(x_i^ax_j^bx_k^c,Z)$. Note that $x_ix_j$ and 
$x_ix_k$ are each covers of the three-cycle, and $x_i^2x_jx_k$ divides $x_i^ax_j^bx_k^c$. 
Let $M_1=x_ix_jN_1$ be any minimal vertex cover of $G$ divisible by $x_ix_j$ but not $x_k$, 
and let $M_2=x_ix_kN_2$ be any minimal vertex cover of $G$ divisible by $x_ix_k$ but not $x_j$. 
Then $x_i^{a-2}x_j^{b-1}x_k^{c-1}M_1M_2 = x_i^ax_j^bx_k^c N_1N_2 \in x_i^ax_j^bx_k^c k[Z] \cap J^2$, 
and therefore $(x_i^ax_j^bx_k^c,Z)$ fails property (b) when $a>1$. A similar argument shows that 
$(x_i^ax_j^bx_k^c,Z)$ fails property (b) if $b>1$ or $c>1$.

Hence $\mult_{J^2}(P_Z) \le 1$ since $(x_ix_jx_k,Z)$ is the only remaining candidate for a standard pair. 
Because $P_Z$ is an associated prime, $\mult_{J^2}(P_Z) \ge 1$, and thus the multiplicity is equal to 1.
\end{proof}  
  
\begin{theorem}\label{adeg}
Let $G$ be a simple graph with $|E_G|$ the number  
of edges of $G$ and $t(G)$ the number of triangles.  
Set $J = I(G)^{\vee}$.
Then $G$ has no odd holes if and only if   
\[\adeg(J^2) = 3|E_G| + t(G).\]  
\end{theorem}  
  
\begin{proof}  
By definition  
\[\adeg(J^2) = \sum_{\mbox{$P \in \Ass(R/J^2)$}} \mult_{J^2}(P).\]  
The associated primes $P$ of $J^2$ are either $P = (x_i,x_j)$ where $x_ix_j \in E_G$,   
$P = (x_i,x_j,x_k)$, where the induced graph on $\{x_i,x_j,x_k\}$ is a triangle (i.e., a 3-cycle),  
or $P = (x_{i_1},\ldots,x_{i_s})$ where the induced graph on  
$\{x_{i_1},\ldots,x_{i_s}\}$ is an odd hole.  Thus  
\begin{eqnarray*}  
\adeg(J^2) &=& \sum_{  
\mbox{$x_ix_j \in E_G$}} \mult_{J^2}((x_i,x_j)) +   
\sum_{\mbox{$\{x_i,x_j,x_k\}$ is a triangle}}   
\mult_{J^2}((x_i,x_j,x_k)) \\  
& & + \sum_{\mbox{$\{x_{i_1},\ldots,x_{i_s}\}$ is an odd hole}} \mult_{J^2}((x_{i_1},\ldots,x_{i_s})).  
\end{eqnarray*}  
So $G$ has no odd hole if and only if   
\begin{eqnarray*}  
\adeg(J^2)&=& \sum_{\mbox{$x_ix_j \in E_G$}} \mult_{J^2}((x_i,x_j)) +   
\sum_{\mbox{$\{x_i,x_j,x_k\}$ is a triangle}}   
\mult_{J^2}((x_i,x_j,x_k)) \\&=& 3|E(G)| + t(G)  
\end{eqnarray*}  
where the final equality follows from Lemma \ref{lemma2} and the calculation just before it.  
\end{proof}  
  
\begin{corollary} \label{cor.deg}  
Let $G$ be a simple graph, and $J = I(G)^{\vee}$.  Then  
\[\deg(J^2) = 3|E(G)|.\]  
\end{corollary}  

\begin{proof}  
The formula for the degree of an ideal $I$ is similar to the 
formula of the arithmetic degree of $I$, except one sums over all minimal associated primes,
instead of all associated primes.  Since the  
minimal associated primes of $J^2$ are precisely those that correspond  
to edges of $G$, the result now follows.    
\end{proof}

\end{document}